\documentclass{article}
\usepackage{amsmath}
\usepackage[amsmath,thmmarks,thref]{ntheorem}
\usepackage[numbers]{natbib}
\include{thmstuff}
\usepackage{amssymb,xspace,amsfonts,latexsym}
\usepackage{graphicx}
\newcommand{\boldp}{\ensuremath{{\bf p}}\xspace}
\newcommand{\boldq}{\ensuremath{{\bf q}}\xspace}
\newcommand{\boldpq}{\ensuremath{{\bf p;\; q}}\xspace}
\newcommand{\boldptimesq}{\ensuremath{{\bf p \times q}}\xspace}
\newcommand{\boldpminq}{\ensuremath{{\bf p; -q}}\xspace}
\newcommand{\chihat}{\ensuremath{\widehat{\chi}}\xspace}
\newcommand{\chitilde}{\ensuremath{\widetilde{\chi}}\xspace}

\newcommand{\compinv}{{\langle -1 \rangle}\xspace}
\newcommand{\ptq}{\ensuremath{p \times q}\xspace}

\begin{document}
\title{Positivity results for Stanley's character polynomials}
\author{A. Rattan\footnote{Current Affiliation:  Department of 
Mathematics, Massachusetts Institute of Technology, email:
{\small\texttt arattan@math.mit.edu}}\\Department of Combinatorics and
Optimization\\University of Waterloo}
\date{May 16, 2006}


\newcommand{\seqtwo}[2]{#1_1, \ldots, #1_{#2}}
\newcommand{\seqinf}[1]{#1_1, #1_{2}, \ldots}
\newcommand{\seq}[2]{\ensuremath{#1_1, #1_2, \ldots, #1_{#2}}}
\newcommand{\seqo}[2]{#1_0, #1_1, \ldots, #1_{#2}}
\newcommand{\seqtwofromto}[3]{#1_{#2}, \ldots, #1_{#3}}
\newcommand{\seqfromto}[3]{#1_{#2}, #1_{#2 + 1}, \ldots, #1_{#3}}
\newcommand{\thbox}{\hfill$\Box$}
\newcommand{\parti}[1]{1^{#1_1}2^{#1_2} \ldots n^{#1_n}}
\newcommand{\parts}[1]{m_1(#1), m_2(#1), \ldots, m_n(#1)}
\newcommand{\partso}[1]{m_0(#1), m_1(#1), \ldots, m_n(#1)}
\newcommand{\set}[2]{\{\, #1 \; |\; #2\}}
\newcommand{\diffpart}[3]{#1_1^{#2_1}#1_2^{#2_2} \ldots #1_{#3}^{#2_{#3}}}
\newcommand{\sumofvars}[2]{#1_1 + #1_2 + \cdots + #1_{#2}}
\newcommand{\sumofvarstwo}[2]{#1_1 + \cdots + #1_{#2}}
\newcommand{\diffofvars}[2]{-#1_1 - #1_2 - \cdots - #1_{#2}}
\newcommand{\lessvars}[2]{#1_1 < #1_2 < \cdots < #1_{#2}}
\newcommand{\greatvars}[2]{#1_1 > #1_2 > \cdots > #1_{#2}}
\newcommand{\prodofvars}[2]{#1_1  #1_2 \cdots #1_{#2}}
\newcommand{\prodtwofromto}[3]{#1_{#2} \cdots #1_{#3}}
\newcommand{\upst}{\mathrm{st}}
\newcommand{\upnd}{\mathrm{nd}}
\newcommand{\upth}{\mathrm{th}}
\newcommand{\coeff}[1]{\ensuremath{[#1]\;}\xspace}

\newcommand{\seqpermvars}[3]{#1_{#2(1)}, #1_{#2(2)}, \ldots, #1_{#2(#3)}}
\newcommand{\smalldisunion}{\ensuremath{\stackrel{\cdot}{\cup}}}
\newcommand{\bigdisunion}{\ensuremath{\stackrel{\cdot}{\bigcup}}}
\newcommand{\symgroup}[1]{\ensuremath{\mathfrak{S}_#1}}
\newcommand{\fullcyc}[1]{\ensuremath{(1 \, 2 \, \cdots \, #1)}}
\newcommand{\laginv}{ {\ensuremath{\langle -1 \rangle}}}
\newcommand{\onenk}{1^{n-k}} 
\newcommand{\supp}{\ensuremath{\mathrm{supp}}\xspace}
\newcommand{\stm}{\ensuremath{\setminus}\xspace}
\newcommand{\al}{\ensuremath{\alpha}\xspace}
\newcommand{\be}{\ensuremath{\beta}\xspace}
\newcommand{\ga}{\ensuremath{\gamma}\xspace}
\newcommand{\lam}{\ensuremath{\lambda}\xspace}
\newcommand{\si}{\ensuremath{\sigma}\xspace}
\newcommand{\om}{{\ensuremath{\omega}\xspace}}
\newcommand{\vt}{\ensuremath{\vartheta}\xspace}
\newcommand{\sms}{\ensuremath{\setminus}\xspace}
\newcommand{\rar}{\ensuremath{\rightarrow}\xspace}
\newcommand{\mpt}{\ensuremath{\mapsto}\xspace}
\newcommand{\f}{\ensuremath{\frac}}
\newcommand{\tf}{\ensuremath{\tfrac}}
\newcommand{\ld}{\ensuremath{\ldots}\xspace}
\newcommand{\mbC}{\ensuremath{\mathbb{C}}\xspace}
\newcommand{\mbR}{\ensuremath{\mathbb{R}}\xspace}
\newcommand{\mcC}{\ensuremath{\mathcal{C}}\xspace}
\newcommand{\mcS}{\ensuremath{\mathcal{S}}\xspace}
\newcommand{\mcR}{\ensuremath{\mathcal{R}}\xspace}
\newcommand{\mcT}{\ensuremath{\mathcal{T}}\xspace}
\newcommand{\mcB}{\ensuremath{\mathcal{B}}\xspace}
\newcommand{\mcP}{\ensuremath{\mathcal{P}}\xspace}
\newcommand{\mfS}{\ensuremath{\mathfrak{S}}\xspace}


\newcommand{\grob}{Gr\"{o}bner basis }
\newcommand{\grobs}{Gr\"{o}bner bases }
\newcommand{\ky}[1]{k[y_1, y_2, \ldots, y_{#1}]}
\newcommand{\kx}[1]{k[x_1, x_2, \ldots, x_{#1}]}
\newcommand{\kyx}[2]{k[y_1, \ldots, y_{#1}, x_1, \ldots, x_{#2}]}
\newcommand{\kxy}[2]{k[x_1, \ldots, x_{#1}, y_1, \ldots, y_{#2}]}
\newcommand{\id}[1]{\langle #1 \rangle}
\newcommand{\qedend}{\hfill\qed}
\newcommand{\isubnk}{I_n^{(k)}}
\newcommand{\psubnk}{P_n^{(k)}}
\newcommand{\khat}{\hat{k}}
\newcommand{\kpark}{{(k)}}
\newcommand{\Phat}{\hat{P}}
\newcommand{\Ihat}{\hat{I}}


\maketitle
\begin{abstract}
\noindent
In \citet{stan:7}, the author introduces expressions for the normalized
characters
of the symmetric group
and states some positivity conjectures for these expressions.  Here, we
give an affirmative partial answer to Stanley's positivity conjectures about the
expressions using results on Kerov polynomials.  In particular, we use new
positivity results in \citet{rattgoul:2}.  We shall see that the generating
series $C(t)$ introduced in \cite{rattgoul:2} is critical to our
discussion.
\end{abstract}

\section{Introduction}\label{sec:intro}

A {\em partition} is a weakly ordered list of positive
integers $\lam =\lam_1\lam_2\ld\lam_k$, where $\lam_1\geq\lam_2\geq\ld\geq\lam_k$.
The integers $\lam_1,\ld ,\lam_k$ are called the {\em parts} of
the partition $\lam$, and we denote the number of parts by $l(\lam)=k$.
If $\lam_1+\ld +\lam_k=d$, then $\lam$ is a partition of $d$, and
we write $\lam\vdash d$. We denote by $\mcP$ the set of all
partitions, including the single partition of $0$ (which
has no  parts).  For partitions $\om, \lam \vdash n$ let $\chi_\om(\lam)$ be
the character of the irreducible representation of the symmetric
group $\mfS_n$ indexed by $\om$, and evaluated on the conjugacy class
 $\mcC_{\lam}$ of $\mfS_n$, which consists of all permutations whose disjoint cycle lengths
are specified by the parts of $\lam$. 

Various scalings of irreducible symmetric group characters have been considered in
the recent literature.  The {\em central character} is given by
\begin{equation*}
\chitilde_{\om}(\lam)=\vert\mcC_{\lam}\vert\f{\chi_\om(\lam)}{{\chi_\om(1^n)}},
\end{equation*}
where $\chi_\om(1^n)$ is
the {\em degree} of the irreducible representation
indexed by $\om$.
For results about the central character, see, for example,~\cite{cgs,fjr,ka:1}.
The scaling to be discussed in this paper, the
                \emph{normalized character}, is given for any
                partitions $\om \vdash n$ and $\mu \vdash k$, where
                $k\leq n$, by
                \begin{equation*}
                        \chihat_\om(\mu\; 1^{n-k}) = n(n-1) \cdots
                        (n-k+1) \frac{\chi_\om(\mu\;
                        1^{n-k})}{\chi_\om(1^n)}.
                \end{equation*}
                For the conjugacy class $C_{k\; 1^{n-k}}$ only, the
                normalized character and the central character are
                related by the following:
                \begin{equation*}
                        \label{eq:chihatdef}
                        \chihat_\om(k1^{n-k}) = n(n-1) \cdots (n-k+1)
                        \frac{ \chi_\om(k1^{n-k})}{\chi_\om(1^n)}=k\chitilde_{\om}(k1^{n-k}).
                \end{equation*}

\begin{figure}[ht]
 \centering
  \includegraphics[width=6.0cm]{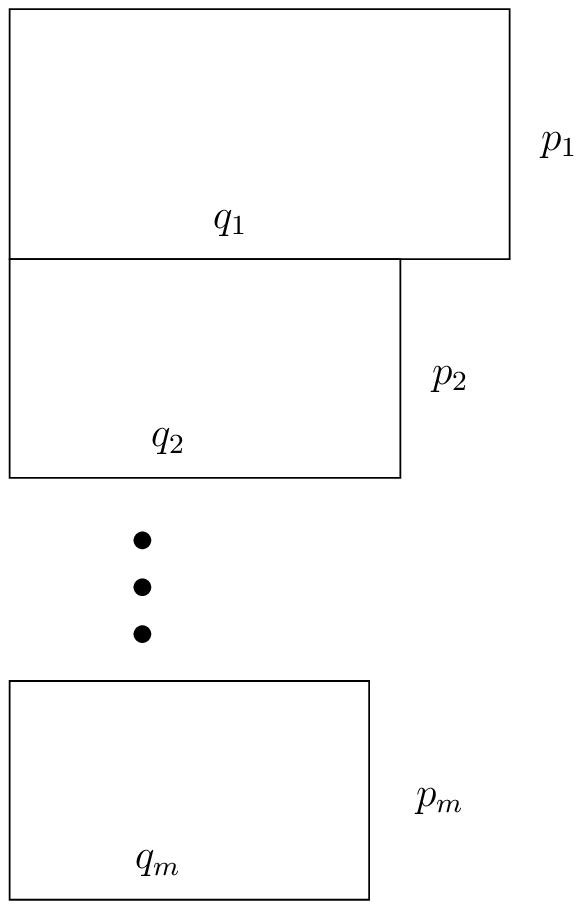}
  \caption{The shape \boldptimesq.}\label{fig:genpart}
\end{figure}

The subject of this paper is a particular polynomial expression for the
normalized character, introduced in \citet{stan:7}.  Consider the 
partition with $p_i$ parts of size $q_i$, for $i$ from 1 to $m$, with $q_1$ the
largest part.
Thus, $\seq{p}{m}$ are
positive integers and $q_1 > q_2 > \cdots > q_m$ (see Figure
\ref{fig:genpart}).  We denote this partition by \boldptimesq.  
Define the series $F_k$ in indeterminates $\seqtwo{p}{m}, \seqtwo{q}{m}$ by
\begin{equation}\label{eq:deffk}
	F_k(\seq{p}{m};\; \seq{q}{m}) = \chihat_\boldptimesq(k\; 1^{n-k})
\end{equation}
We often denote
$(\seqtwo{p}{m})$ by $\boldp$ and $(\seqtwo{q}{m})$ by \boldq, giving us the
notation $F_k(\boldp;\; \boldq)$
for $F_k(\seq{p}{m};\; \seq{q}{m})$.  The following theorem appears in \citet[Proposition 1]{stan:7}.
\begin{theorem}[Stanley]
$F_k(\boldp;\; \boldq)$ is a polynomial in the $p$'s and $q$'s such that
\linebreak $F_k(1,1,\dots,1;\; -1,-1,\dots,-1) = (k+m-1)_k$. 
\end{theorem}
In light of this theorem, we call the polynomials in \eqref{eq:deffk}
\emph{Stanley's character polynomials}.  These polynomials are the main objects
in this paper.  For example, for the case $m=2$, the first two Stanley
polynomials are
\begin{align*}
F_1(a,p;\; b, q) &= -ab-pq,\\
F_2(a,p;\; b, q) &= -{a}^{2}b+a{b}^{2}-2\,apq-{p}^{2}q+p{q}^{2}
\end{align*}
where we have set $p_1 = a, p_2 = p, q_1 = b$ and $q_2 =q$.

\subsection{Main results} \label{sec:intromain}

In \cite{stan:7}, Stanley generalizes $F_k(\boldp;\; \boldq)$ to
\begin{equation*}
	F_\mu(\boldp;\; \boldq) = \chihat_{\boldptimesq}(\mu\;
	1^{n-k}),
\end{equation*} 
where $\mu$ is a partition of $k$.  Stanley states that $F_\mu(\boldp;\; \boldq)$ is, by the Murnaghan-Nakayama rule, a polynomial with integer coefficients.
In \cite[Conjecture 1]{stan:7}, Stanley gives a
positivity conjecture for a variant of the series $F_\mu(\boldp;\; \boldq)$.  For
convenience, we use the notation $-\boldq = (-q_1, -q_2, \dots, -q_m)$, and
$F_\mu(\boldp;\; -\boldq)$ is the series $F_\mu(\boldp;\; \boldq)$ with $q_i$
replaced by $-q_i$.  Call any series $T(\boldp;\; \boldq)$ in the indeterminates
$p$'s and $q$'s \emph{$\boldp, \boldq$-positive} if the coefficients of all
terms are positive.
\begin{conjecture}[Stanley] \label{conj:main}
	For any partition $\mu \vdash k$, the series $(-1)^kF_\mu(\boldpminq)$
	is $\boldp, \boldq$-positive.
\end{conjecture}
Stanley only proves this in the case $m=1$, the so-called rectangular case as in
this case the shape $\boldptimesq$ is the rectangle with $p_1$ parts all equal
to $q_1$.  We drop the subscript 1 in this case and say $p \times q$ is the
partition with $p$ parts all equal to $q$.  In the rectangular case, Stanley
proves positivity by giving a stronger result; he gives a combinatorial
interpretation for the coefficients, given in
\cite[Theorem 1]{stan:7} and stated below.
\begin{theorem}[Stanley]
	\label{thm:stantwovar}  Suppose that $\ptq \vdash n$ and $\mu \vdash k$
	for $k \leq n$.  Let $\lam_\mu$ be any fixed permutation in the
	conjugacy class indexed by $\mu$ in $\mfS_k$.  Then,
	\begin{equation*}
		\chihat_{\ptq}(\mu\; 1^{n-k}) = (-1)^k \sum_{u, \nu \atop u \nu =
		\lam_\mu} p^{\ell(u)} (-q)^{\ell(\nu)}.
	\end{equation*}
\end{theorem}
For general $m$, Conjecture \ref{conj:main} remains open.  In fact, there is no
proof even in the case where $\mu$ has one part, that is, it is not yet known
whether $(-1)^k F_k(\boldpminq)$ is $\boldp, \boldq$-positive.  This generating
series will be the focus of this paper and we address its $\boldp,
\boldq$-positivity.

Stanley does state that the terms of highest degree, the terms of degree $k+1$, of
$(-1)^kF_k(\boldpq)$ have a particularly nice form, and are given in
\eqref{eq:stantop} below.  He does not, however, prove that these terms are
$\boldp, \boldq$-positive but does state that Elizalde has given a proof of this
in private communication to him (see \eqref{eq:elizalde} below).  The proof by
Elizalde does not appear to be anywhere in the literature.

In this paper, we give a new proof of Theorem \ref{thm:stantwovar}.  We do this
using \emph{shift symmetric functions}.  This new proof, we hope, is simpler and
makes the result more transparent.  Furthermore, it highlights the already known
connection between shift symmetric functions and the normalized character.  As
for the general case of $F_k(\boldpq)$, 
we give a proof of $\boldp,\boldq$-positivity of the terms of highest degree,
the terms of degree $k+1$,
in $(-1)^k F_k(\boldpminq)$ (as mentioned above, this was also proved by
Elizalde, but his proof does not appear in the literature).  We also give a proof that the terms of degree
$k-1$ and $k-3$ in $(-1)^k F_k(\boldpminq)$, the terms of second and third
highest degree, are also $\boldp, \boldq$-positive, which are new results.  We do this by using recent
results concerning \emph{Kerov's polynomials}.  For Kerov's polynomials,
there is a notion of R-positivity and a new notion of C-positivity introduced in
\citet{rattgoul:2}, which we shall use to show our positivity results for
$(-1)^k F_k(\boldpminq)$.  Finally, we end the paper by showing 
that C-positivity of Kerov's polynomials implies $\boldp,
\boldq$-positivity of $(-1)^k F_k(\boldpminq)$, also a new result.

The necessary results on Kerov's polynomials are reviewed in the next section.

\subsection{Kerov polynomials}
We adapt the following description from
Biane~\cite{bi1,bi2}: consider the Young diagram of $\om$, in the French
convention (see \cite[footnote page 2]{mac:1}),
and translate it, if necessary, so that the
bottom left of the diagram is placed at the origin of an $(x,y)$ plane.
Finally, 
rotate the diagram counter-clockwise by $45^\circ$.  Note that $\om$ is uniquely
\begin{figure}[ht]\label{fig:secondpart}
 \centering
  \includegraphics[width=12.0cm]{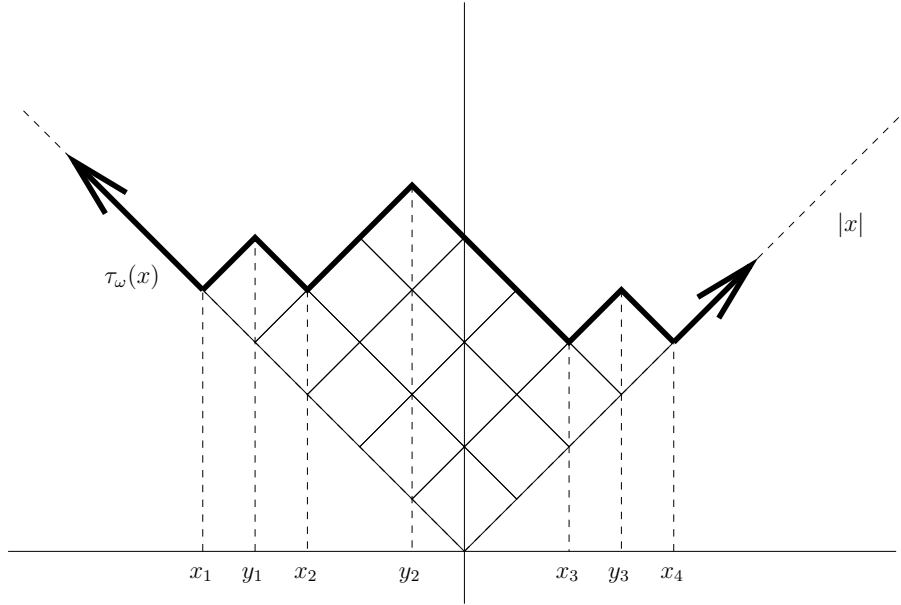}
  \caption{The partition (4\,3\,3\,3\,1) of 14, drawn in the French convention,
  and rotated by $45^\circ$.}
\end{figure}
determined by the curve $\tau_\om(x)$ (see
Figure \ref{fig:secondpart}).  The value of
$\tau_\om(x)$ is equal to $|x|$ for large negative or positive values of $x$ and it
is clear that $\tau_\om^\prime (x) = \pm 1$, where differentiable. The
points $x_i$ and $y_i$ are the
$x$-coordinates of the local minima and maxima, respectively, of the
curve $\tau_\om(x)$.  We suitably scale the size of the boxes in our Young
diagram so that the points $x_i$ and $y_i$ are integers.
Setting $\sigma_\om(x) = (\tau_\om(x) - |x|)/2$, consider the function
\begin{equation*}
  H_\om(z) =\f{1}{z} \exp \int_\mathbb{R} \frac{1}{z-x}
  \; \sigma_\om^\prime(x) \; dx .
\end{equation*}
Carrying out the above integration one obtains
\begin{equation}\label{eq:seriesh}
H_\om(z) = \frac{\prod_{i=1}^{m-1} (z - y_i)}{\prod_{i=1}^m (z - x_i)} ,
\end{equation}
where $m$ is the number of nonempty rows in the Young diagram of $\om$.
Now let $R_\om(z) = \sum_{i\geq 0} R_i(\om) z^i$ be defined by 
\begin{equation}
R_\om(z) = z H_{\om}^{\langle -1
\rangle}(z),\label{eq:defofr}
\end{equation}
where $\langle -1\rangle$ denotes compositional inverse.
We will be applying Lagrange inversion (see \citet[Section 1.2]{goul:1} or
\citet[Theorem 5.4.2]{sta:2}) to
\eqref{eq:seriesh} and \eqref{eq:defofr}, in which case we change
$H_\om(z)$ to a formal power series.  We then obtain
\begin{equation}\label{eq:rinh}
		R_\om(z) =
		\frac{z}{\left(H_\om(1/z)\right)^\compinv}.
\end{equation}
Briefly, the origins of the series $H_\om(z)$, and in fact 
Kerov's polynomials, come from
attempting to answer asymptotic questions about the characters of the
symmetric group using free probability.  In that context $H_\om(z)$ is called
the \emph{moment generating series} (traditionally denoted $G_\om(z)$) and
$R_\om(z)$ is the \emph{free cumulant generating series} (traditionally denoted
$K_\om(z)$).  We refer the reader to \citet{bi, bi2} for the background on
asymptotics of characters and free probability.

Finally, the polynomials we are concerned with involve the $R_i(\om)$'s and
are given in the following theorem.  They first appeared in \citet[Theorem
1.1]{bi2}.
\begin{theorem}[Biane]\label{kebi}
For $k\geq 1$, there exist universal polynomials $\Sigma_k$, with
  integer coefficients, such that
  \begin{equation}
    \label{eq:kerovpolys}
    \chihat_\omega(k\; 1^{n-k}) = \Sigma_k(R_2(\omega), R_3(\omega), \ldots,
    R_{k+1}(\om)),
  \end{equation}
for all $\om\vdash n$ with $n\geq k$.
\end{theorem}
Biane attributes Theorem~\ref{kebi} to Kerov, who described this result
in a talk at an IHP conference in 2000, but a proof first appears
in a later paper of Biane~\cite{bi}. The polynomials $\Sigma_k$
are known as {\em Kerov's character polynomials}.
They are referred to as ``universal polynomials'' in Theorem~\ref{kebi}
to emphasize that they are independent of $\om$ and $n$, subject
only to $n\geq k$. Thus we write them with $R_i(\om)$ replaced by
an indeterminate $R_i$, $i\geq 2$.
In indeterminates $R_i$,
the first six of Kerov's character polynomials, as listed in~\cite{bi2}, are
given below:
\begin{align}
  \Sigma_1 & = R_2\nonumber\\
  \Sigma_2 & = R_3\nonumber\\
  \Sigma_3 & = R_4 + R_2\nonumber\\
  \Sigma_4 & = R_5 + 5R_3\label{eq:rexp}\\
  \Sigma_5 & = R_6 + 15R_4 + 5R_2^2 + 8R_2\nonumber\\
  \Sigma_6 & = R_7 + 35R_5 + 35R_3 R_2 + 84R_3\nonumber
\end{align}
We note that once one knows Kerov's polynomial $\Sigma_k$ in order to find
$\chihat_\om(k\; 1^{n-k})$ all one needs is to construct the series $H_\om(z)$
as in \eqref{eq:seriesh}
and apply Lagrange inversion to \eqref{eq:rinh} to find the series $R_\om(z)$, giving the $R_i(\om)$.
Substituting these values into the $k^\upth$ polynomial will then give the
normalized character $\chihat_\om(k\; 1^{n-k})$.

Note that all coefficients appearing in this list are positive.
It is conjectured that this holds in general: that for any $k\geq 1$, all nonzero
coefficients in $\Sigma_k$  are positive (see Theorems \ref{thm:kerovtopterm},
\ref{thm:sigk2} and \ref{thm:cpos} for positivity results obtained so far for
Kerov's polynomials).  Kerov's polynomials remain somewhat
of a mystery, in spite of recent efforts by \citet{bi, sniasymp} and \citet{rattgoul:2}.  In
particular, in \cite{rattgoul:2} the authors introduce the
polynomial 
$C(t)=\sum_{m\geq 0}C_mt^m$ given by 
\begin{equation}\label{eq:Cdef}
C(t)=\f{1}{1-\sum_{i\geq 2}(i-1)R_it^i}.
\end{equation}
From the definition of $C(t)$ we see the $C_m$ are
polynomials in the $R_i$'s, with $C_0=1$, $C_1=0$, and
\begin{equation}\label{cexplicit}
C_m=\sum_{\stackrel{j_2,j_3,\ld \geq 0}{2j_2+3j_3+\ld =m}}
\!\!\!\!\!\!\!\! (j_2+j_3+\ld )!\prod_{i\geq 2}\f{((i-1)R_i)^{j_i}}{j_i!},
\;\;\;\; m\geq 2.
\end{equation}
Writing Kerov's polynomials in terms of the $C$'s we have (from
\citet[page 7]{rattgoul:2})
\begin{align}
\Sigma_3-R_4&=C_2\nonumber\\
\Sigma_4-R_5&=\tf{5}{2}C_3\nonumber\\
\Sigma_5-R_6&=5\,{C}_{{4}}+8\,{C}_{{2}}\label{eq:cexp}\\
\Sigma_6-R_7&=\tf{35}{4}\,{C}_{{5}}+42\,{C}_{{3}}\nonumber\\
\Sigma_7-R_8&= 14\,{C}_{{6}}+{\tf{469}{3}}\,{C}_{{4}}
+{\tf{203}{3}}\,{{C}_{{2}}}^{2} + 180\,{C}_{{2}}\nonumber\\
\Sigma_8-R_9&=21\,{C}_{{7}}+{\tf{1869}{4}}\,{C}_{{5}}
+{\tf{819}{2}}\,C_3 C_2+1522\,{C}_{{3}}\nonumber
\end{align}
We will, henceforth, call the expansions in \eqref{eq:rexp} and \eqref{eq:cexp}
the \emph{R-expansions} and \emph{C-expansions}, respectively, of Kerov's
polynomials.  We will also call the property that all coefficients of the
$R$'s are positive \emph{R-positivity} and \emph{C-positivity} is analogously
defined for $\Sigma_k - R_{k+1}$.  It follows from \eqref{cexplicit} that
C-positivity of Kerov's polynomials implies R-positivity.

Define the weight of a monomial $R_{i_1}^{j_1} R_{i_2}^{j_2} \cdots
R_{i_m}^{j_m}$ to be $\sum_{t = 1}^m i_t j_t$ (and analogously for monomials in
$C$'s).  Further let $\Sigma_{k, 2n}$ be the terms of weight $k + 1 - 2n$ in
$\Sigma_k$.  From the combinatorial origins of Kerov's polynomials, it follows
that in the $k^\upth$ Kerov polynomial $\Sigma_k$ the terms of weight $k\,
(\mathrm{mod}\; 2)$ each have zero coefficient (see \citet[Section 4]{bi}).  The
following theorems appear in \citet[Proof of Theorem 1.1]{bi}, \citet[Theorems
1.3 and 3.3]{rattgoul:2} and \citet[Section 1.3.3]{sniasymp}.
\begin{theorem}[Biane]\label{thm:kerovtopterm}
	\begin{equation*}
		\Sigma_{k,0} = R_{k+1}.
	\end{equation*}
	That is, there is only one term of weight $k+1$ in $\Sigma_{k}$ and it
	is $R_{k+1}$.
\end{theorem}
\begin{theorem}[Goulden-R, \'Sniady]\label{thm:sigk2}
	\begin{equation*}
		\Sigma_{k,2} = \f{1}{4}{k+1 \choose 3} C_{k-1}.
	\end{equation*}
	Consequently, $\Sigma_{k,2}$ is C-positive.
\end{theorem}
\begin{theorem}[Goulden-R]\label{thm:cpos}
	$\Sigma_{k,4}$ is C-positive.	
\end{theorem}
Finally, the following theorem is a corollary of \cite[Theorem 2.1]{rattgoul:2}.

\begin{theorem}\label{thm:gencexp}
	For $k \geq 1$,
	\begin{equation*}
		\Sigma_{k,2n} = \sum_{i_1, i_2,\dots, i_{2n-1} \geq 0 \atop i_1
		+ i_2 + \cdots + i_{2n-1} = k+1 -2n} \gamma_{i_1,
		i_2, \dots, i_{2n-1}} C_{i_1} \cdot C_{i_2} \cdots
		C_{i_{2n-1}}
	\end{equation*}
	where the $C_t$ are given in \eqref{eq:Cdef} and the $\gamma$'s are
	rational.  In particular,
	$\Sigma_{k,2n}$
	is C-positive (and, consequently, R-positive)
	\index{C-positivity}\index{R-positivity} if all $\gamma_{i_1,
		i_2, \dots, i_{2n-1}}$ are positive.
\end{theorem}

\section{Stanley's Polynomials for Rectangular Shapes}\label{sec:rect}

In this section we study a specific two variable case
of Stanley's results as they have a particularly beautiful form. 

We begin with the normalized character $\chihat_\om$ when $\om$ has the rectangular shape
of $p$ parts, all equal to $q$.  In Section \ref{sec:intromain}, we denoted this
shape with $\ptq$.  Further, in Section \ref{sec:intromain} we gave the central
theorem in this case as Theorem \ref{thm:stantwovar}.  
This result can be written in terms of the \emph{connection coefficients} of the
symmetric group;  Theorem \ref{thm:stantwovar} then becomes
\begin{equation*}
	\chihat_{\ptq}(\mu\; 1^{n-k}) = (-1)^k \sum_{u, \nu \vdash k}
		c_{u, \nu}^\mu p^{\ell(u)} (-q)^{\ell(\nu)}.
\end{equation*}
Here the $c^\mu_{u, \nu}$ are defined as the structure constants of
the central elements $K_u$ of the group algebra of $\mfS_n$;  that is,
\begin{equation*}
	K_u K_\nu = \sum_{\mu} c^\mu_{u, \nu} K_\mu.
\end{equation*}
Stanley's proof of this involves a combination of results;  Stanley uses results about certain tableaux, the Murnaghan-Nakayama rule, and the following symmetric function identity
\begin{equation*}
	\sum_{\om \vdash k} H_\om\; s_\om(x) s_\om(y) s_\om(z) = \sum_{\om \vdash
	k} p_\om(x) p_\om(y) p_\om(z),
\end{equation*}
which appears in \citet{han:1} (here, $s_\om(x)$ and $p_\om(x)$ are the Schur
symmetric function and power sum symmetric function, respectively).
Here, we present an original proof with the aim of making the result
more transparent and, in addition, of showing more connections between what are
known as \emph{shift symmetric functions} and the normalized character
$\chihat_{\ptq}$
(we shall see that there is already a known relationship between these
objects).  Sections \ref{sec:shiftsym} gives the necessary background on shift
symmetric functions for this proof.

\subsection{A Brief Account of Shift Symmetric Functions}\label{sec:shiftsym}

In \citet{ok:1}, the authors define \emph{shift Schur polynomials} as
\begin{equation*}
	s_\lam^*(\seq{x}{n}) = \f{\det\left((x_i + n -i)_{\lam_j + n
	-j}\right)_{1 \leq i,j \leq n}}{\det\left((x_i)_{n-j}\right)_{1 \leq i,j
	\leq n}}.
\end{equation*}
The \emph{shift Schur functions}
denoted by $s_\lam^* \in \Lambda^*$
are defined as the inverse limit of the sequence $(s_\lam^*(\seq{x}{n}))_{n \geq
1}$.   Further, one can define the \emph{p-sharp} shift symmetric functions $p_\mu^\sharp$;  they are
\begin{equation*}
	p_\mu^\sharp = \sum_{\rho \vdash k} \chi_\rho(\mu) s_\rho^*
\end{equation*}
(see \citet[Section 1]{ok:1} for more details).

The following result connects shift symmetric functions and the scaled
characters $\chihat$, and can be found in \citet[(15.21)]{ok:1}.
\begin{theorem}[Okounkov, Olshanski]
	\label{thm:psharp}  Suppose that $\mu \vdash k$ and $\lam \vdash n$.
	Then
	\begin{equation*}
		p_\mu^\sharp(\lam) = \chihat_\lam(\mu\; 1^{n-k}),
	\end{equation*}
	where $p^\sharp(\lambda)$ is the substitution $x_i = \lam_i$ for $1 \leq
	i \leq \ell(\lam)$ and $x_i = 0$ for $i > \ell(\lam)$.
\end{theorem}
The following theorem gives a combinatorial interpretation to shift Schur
functions; it is also found in \citet[Theorem 11.1]{ok:1}.  For any shape $\mu$,
a \emph{reverse tableau of shape $\mu$} is
a function $T : \textnormal{boxes of } \mu
\mapsto \mathbb{P}$, where $\mathbb{P}$ is the set of positive integers,
such that $T$ is weakly decreasing along the rows of $\mu$ and strongly
decreasing along the columns of $\mu$.  We denote by RTab$(\mu)$ the set of
reverse tableau of shape $\mu$.
\begin{theorem}[Okounkov, Olshanski]
	\label{thm:combshiftschur} For $\lambda \in \mathcal{P}$, 
	\begin{equation*}
		s_\lam^* = \sum_{T \in \mathrm{RTab}(\mu)} \prod_{u \in \mu}
		(x_{T(u)} - c(u)),
	\end{equation*} where $T(u)$ is the value assigned to the box $u$ by the
	tableau $T$ and, again, $c(u)$ is the content of the box $u$.
\end{theorem}

\subsection{Proof of Theorem \ref{thm:stantwovar}}\label{sec:origproofstan}

We are now ready to give a proof of Theorem \ref{thm:stantwovar}.

\noindent
{\bf Proof of Theorem \ref{thm:stantwovar}.}
As a first step to this proof, for a partition $\lambda \vdash k$ we evaluate $s_\lambda^*(\seq{x}{p})$ with $x_i =
q$ for $1 \leq i \leq p$;  that is, we compute the evaluation
$s_\lambda^*(\ptq)$.  
Using Theorem \ref{thm:combshiftschur} we obtain
\begin{align}
	s_\lambda^*(\ptq) &= \left. \sum_{T \in \mathrm{RTab}(\lambda)} \prod_{u
	\in \lambda} (x_{T(u)} - c(u)) \right|_{(\seqtwo{x}{p}) =
	(q,\dots,q)}\nonumber\\
	&= \sum_{T \in \mathrm{RTab}(\lambda)} \prod_{u \in \lambda} (q -
	c(u))\nonumber\\
	&= \left.(-1)^k \prod_{u \in \lambda} (-q + c(u)) \sum_{T \in
	\mathrm{RTab}(\lambda)} 1\right|_{(\seqtwo{x}{p}) =
	(q,\dots,q)}.\label{eq:staralmost} 
\end{align}
The number of $\mathrm{RTab}(\lambda)$ is clearly the number of
\emph{semi-standard Young tableaux} (see \citet[page 309]{sta:2}) of
shape $\lambda$ filled with only numbers $1,2,\dots, p$, which is 
$s_\lambda({\bf 1^p})$ from the combinatorial definition of Schur functions.  Thus, from
\eqref{eq:staralmost} above and the well known specialization of the Schur
functions
\begin{equation*}
	s_\lambda({\bf 1^p}) = \f{\prod_{u \in \lambda} (p +
	c(u))}{H_\lambda},
\end{equation*}
where $s_\lambda({\bf 1^p})$ is obtained by setting $x_i = 1$ for all $1 \leq i
\leq p$ and $x_i = 0$ for all $i > p$ in the Schur function $s_\lambda({\bf x})$.
\begin{align*}
	s_\lambda^*(\ptq) &= (-1)^k\; \prod_{u \in \lambda} (-q + c(u))\;
	s_\lambda({\bf 1^p})\\
	&= \f{(-1)^k}{H_\lambda} \; \prod_{u \in \lambda} (-q + c(u)) (p +
	c(u)).
\end{align*}
Therefore, from Theorem \ref{thm:psharp} and \eqref{eq:staralmost} we have
\begin{align}
	\chihat_{\ptq}(\mu\; 1^{n-k}) &= \sum_{\lambda \vdash k}
	\chi_{\lambda}(\mu) s_\lambda^*(\ptq)\nonumber\\
	&= (-1)^k \sum_{\lambda \vdash k} \f{\chi_\lambda(\mu)}{H_\lambda}
	\prod_{u \in \lambda} (p + c(u))(-q + c(u)) \nonumber\\
	&= (-1)^k \sum_{\alpha, \beta, \lambda \vdash k}
	\f{\chi_{\lambda}(\mu)}{H_\lambda}
	\f{|C_\alpha|}{f^\lam} \chi_\lam(\alpha) p^{\ell(\alpha)}
	\f{|C_\beta|}{f^\lam} \chi_\lam(\beta) (-q)^{\ell(\beta)}
	\nonumber\\
	&= (-1)^k \sum_{\alpha, \beta, \vdash k}
	p^{\ell(\alpha)} (-q)^{\ell(\beta)}
	\f{|C_\alpha||C_\beta|}{k!} \sum_{\lambda \vdash k}
	\f{1}{f^\lam} \chi_\lam(\alpha) \chi_\lam(\beta)\chi_{\lambda}(\mu)
	\nonumber\\
	&= (-1)^k \sum_{\alpha, \beta \vdash k} p^{\ell(\alpha)}
	(-q)^{\ell(\beta)} c^\mu_{\alpha,
	\beta},\nonumber
\end{align}
where the third equality follows from the well known identity
\begin{equation}\label{eq:contentprod}
		\prod_{u \in \lam} (x + c(u)) = \sum_{\beta \vdash n}
		\f{|C_\beta|}{f^\lam} \chi_\lam(\beta) x^{\ell(\beta)}.
\end{equation}
and the last equality follows from the well known identity
\begin{equation*}
	\phantom{} [K_\mu] K_\alpha K_\beta = \f{|C_\alpha| |C_\beta|}{k!}
	\sum_{\lam \vdash k} \f{1}{f^\lam} \chi_\lam(\alpha)
	\chi_\lam(\beta)\chi_\lam(\mu)
\end{equation*}
(see \citet[Lemma 2.4]{jack:1}). This completes the proof.\hfill$\Box$

\section{Generalizations to Non-Rectangular Shapes}\label{sec:gener}

We now deal with the case of a general shape $\boldptimesq$ but, as mentioned in
Section \ref{sec:intromain}, we are concerned with the series $F_\mu(\boldpq)$
when $\mu$ has a single part;  that is, we are only concerned with the series
$F_k(\boldpq)$.  
The expressions $(-1)^kF_k(\boldp;\; -\boldq)$ for $k = 1,2,3,4$ and $m=2$
are given in \eqref{eq:exstanpolys}.  These data also appear in \citet[page
8]{stan:7}.
\begin{align}
-F_1(a,p;\; -b, -q) &= ab+pq,\nonumber\\
F_2(a,p;\; -b, -q) &= {a}^{2}b+a{b}^{2}+2\,apq+{p}^{2}q+p{q}^{2},\nonumber\\
-F_3(a,p;\; -b, -q) &= 
{a}^{3}b+3\,{a}^{2}{b}^{2}+3\,{a}^{2}pq+a{b}^{3}+3\,abpq
+3\,a{p}^{2}q\nonumber\\
& \;\;\;\; +3\,ap{q}^{2}+{p}^{3}q+3\,{p}^{2}{q}^{2}+p{q}^{3}+ab+pq,\nonumber\\
F_4(a,p;\; -b, -q) &= 
{a}^{4}b+6\,{a}^{3}{b}^{2}+4\,{a}^{3}pq+6\,{a}^{2}{b}^{3}+12\,{a}^{2}b
pq\label{eq:exstanpolys}\\
& \;\;\;\; +6\,{a}^{2}{p}^{2}q+6\,{a}^{2}p{q}^{2}+a{b}^{4}+4\,a{b}^{
2}pq+4\,ab{p}^{2}q\nonumber\\
& \;\;\;\; +4\,abp{q}^{2}+4\,a{p}^{3}q+14\,a{p}^{2}{
q}^{2}+4\,ap{q}^{3}+{p}^{4}q\nonumber\\
& \;\;\;\; +6\,{p}^{3}{q}^{2}+6\,{p}^{2}{q}^{3
}+p{q}^{4}+5\,{a}^{2}b+5\,a{b}^{2}+10\,apq+5\,{p}^{2}q\nonumber\\
& \;\;\;\; +5\,p{q}^{2}.\nonumber
\end{align}
Stanley mentions that the terms of highest degree in $F_k(\boldp;\;
\boldq)$, \emph{i.e.} the terms of degree $k+1$, have a particularly nice expression.
Keeping with Stanley's notation, let $G_k(\boldp;\;
\boldq)$
be the terms of
highest degree in $F_k(\boldp;\; \boldq)$.  We have the following expression for
the generating series of $G_k(\boldp;\; \boldq)$, which we 
call $G_{\boldpq}(z)$
This theorem appears, with proof, in \cite[Proposition 2]{stan:7}.
\begin{theorem}[Stanley]  The generating series for $G_k(\boldp;\; \boldq)$ is
\begin{equation}
	G_{\boldpq}(z) = 1 + \sum_{i \geq 1} G_{i-1}(\boldp;\; \boldq)
	z^i 
	= \f{z}{\left(\displaystyle\f{\displaystyle z
	\prod_{i=1}^m\left(1-\left(q_i + \sum_{j=i+1}^m
	p_j\right)z\right)}{\displaystyle\prod_{i=1}^m \left(1 - \left(q_i +
	\sum_{j=i}^m p_j \right)z\right)} \right)^\laginv}\label{eq:stantop}.
\end{equation}
\end{theorem}
Of course, $\boldp, \boldq$-positivity
of $(-1)^k F_k(\boldp;\; -\boldq)$ would imply that
$(-1)^k G_k(\boldp;\; -\boldq)$ is also $\boldp, \boldq$-positive.  Stanley does
not prove $\boldp, \boldq$-positivity for the latter series in \cite{stan:7} but states that S. Elizalde has proven this in a
private communication to him.  In fact, Elizalde shows (according to
Stanley)
\begin{multline}
	(-1)^k G_k(\boldp;\; \boldq) = \f{1}{k} \sum_{\sumofvarstwo{i}{m} +
	\sumofvarstwo{j}{m} = k+1} {k \choose i_1} \left( {i_1 \choose j_1}
	\right)\\
	\prod_{s=2}^m \left( \sum_{r=0}^{\mathrm{min}(i_s,j_s)} {k
	\choose r} \left( {r \choose j_s -r} \right) {k - r - i_i - \cdots -
	i_{s-1} - j_1 - \cdots - j_{s-1} \choose i_s -r}
	\right) \label{eq:elizalde} \\
	\cdot p_1^{i_1} \cdots
	p_m^{i_m} q_1^{i_1} \cdots q_m^{i_m},
\end{multline}
where $\left( {n \choose k} \right) = {n + k -1 \choose k}$.
However, as far as this author can see, no proof
exists in the literature.  

In the next sections we give partial answers to the positivity
questions concerning $(-1)^k F(\boldp;\; -\boldq)$.  As alluded to 
in the Section \ref{sec:intro}, we use Kerov's polynomials to answer these
questions.

\section{Applying Kerov Polynomials to Stanley's Polynomials}\label{sec:kerovstan}

Note that both Kerov's polynomials, along with \eqref{eq:defofr}, and \eqref{eq:deffk} give expressions
for the scaled character $\chihat_\om$.  Since they hold for any
shapes $\boldptimesq$, we can conclude that they give the same expression for
$\chihat_\om$.  Thus, we will use \eqref{eq:defofr} and \eqref{eq:kerovpolys} 
to obtain results about Stanley's polynomials.  More specifically, using
\eqref{eq:defofr} we obtain the $R_i$ in Kerov's polynomials for the general
shape $\boldptimesq$;  we then use the $R_i$ along with Theorems
\ref{thm:kerovtopterm}, \ref{thm:sigk2} and \ref{thm:cpos} to give some
positivity results for Stanley's polynomials.  
The main theorem needed to give our positivity results is given in Theorem
\ref{thm:phiform} of Section \ref{sec:kplusoneterms};  also in Section
\ref{sec:kplusoneterms} we show, using Theorems \ref{thm:phiform} and
\ref{thm:kerovtopterm}, that the terms of highest degree of Stanley's polynomials
are positive.  In Section \ref{sec:kminoneterms} we use Theorem
\ref{thm:phiform} and Theorems \ref{thm:sigk2} and \ref{thm:cpos} to prove
the positivity of the terms of degree $k-1$ and $k-3$ in $F_k(\boldp;\;
\boldq)$.  Finally, we end the paper by showing in Theorem
\ref{thm:connectkerovstan} that C-positivity for Kerov's
polynomials implies $\boldp, \boldq$-positivity for Stanley's polynomials.

\subsection{The Series $H$ for the Shape \boldptimesq}

We now compute what the series $H$ in \eqref{eq:seriesh} must be for the shape \boldptimesq.  For the shape
\boldptimesq, it is not difficult to see that its interlacing sequence of maxima and
minima is
\begin{equation*}
	x_1 = q_1,\;\;\; y_1 = q_1 - p_1,\;\;\; x_2 = q_2 - p_1,\;\;\; y_2 =
	q_2 - p_1 - p_2,\;\;\; x_3 = q_3 - p_1 - p_2,
\end{equation*}
\begin{equation*}
	y_3 = q_3 -p_1 - p_2 - p_3,\;\;\; \dots,
	\;\;\; x_{m-1} = q_m - \sum_{i=1}^{m-1} p_i,\;\;\; y_m = q_m -
	\sum_{i=1}^m p_i,\;\;\; x_m = -\sum_{i=1}^m p_i.
\end{equation*}
From \eqref{eq:seriesh}, we have
\begin{align}
	H_{\boldptimesq}(1/z) &= \f{\displaystyle z \left(1 - \left(q_1 - p_1\right)z\right) \left(1 - \left(q_2 - \left(p_1 + p_2\right)\right)z\right) \cdots
	\left(1 - \left(q_m - \sum_{i=1}^m
	p_i\right)z\right)}{\displaystyle \left(1 - q_1 z\right) \left(1 - \left(q_2 - p_1\right) z\right) \cdots
	\left(1 - \left(q_m - \sum_{i=1}^{m-1} p_i\right)z\right)\left(1 + \sum_{i=1}^m p_i\right)}\nonumber\\
	&= \f{\displaystyle z \prod_{i=1}^m \left(1 - \left(q_i - \sum_{j=1}^i
	p_j\right)z\right)}{\displaystyle \left(1 + \sum_{j=1}^m
	p_j z\right) \prod_{i=1}^m \left(1-\left(q_i - \sum_{j=1}^{i-1} p_j\right)
	z\right)},\label{eq:stanphi}
\end{align}
and we obtain from \eqref{eq:rinh}
\begin{equation}\label{eq:defstantopkerov}
	R_{\boldptimesq}(z) = \f{z}{\left(\f{\displaystyle z \prod_{i=1}^m \left(1 -
	\left(q_i -
	\sum_{j=1}^i p_j\right)z\right)}{\displaystyle \left(1 + \sum_{j=1}^m
	p_j z\right) \prod_{i=1}^m \left(1-\left(q_i - \sum_{j=1}^{i-1}
	p_j\right) z \right)} \right)^\laginv}.
\end{equation}
Alternatively, it follows from Lagrange inversion (see
\citet[Section 2]{goul:1} or \citet[Theorem 5.4.2]{sta:2}) that if
\begin{align}
	\phi_\boldptimesq(z) &= \f{z}{H_{\boldptimesq}(1/z)}\nonumber\\
	&= \f{\left(1 + \sum_{j=1}^m
	p_j z\right) \prod_{i=1}^m \left(1-\left(q_i - \sum_{j=1}^{i-1} p_j\right) z\right)}{\prod_{i=1}^m
	\left(1 - \left(q_i - \sum_{j=1}^i p_j\right)z\right)},\label{eq:phipq}
\end{align}
then
\begin{equation}\label{eq:recr}
	\f{z}{R_\boldptimesq(z)} = z\, \phi_\boldptimesq
	\left(\f{z}{R_\boldptimesq(z)}\right).
\end{equation}
Applying Lagrange inversion, we obtain for $k \geq 2$
\begin{align}
	R_k(\boldptimesq) &= [z^{k-1}] \f{R(z)}{z}\nonumber\\
	&= \f{1}{k-1} [y^{k-2}] -\f{1}{y^2} \phi_{\boldptimesq}^{k-1}(y)\nonumber\\
	&= -\f{1}{k-1} [y^k] \phi_{\boldptimesq}^{k-1}(y). \label{eq:altrk}
\end{align}
Of course, substituting $R_i(\boldptimesq)$ for $R_i$ in Kerov's polynomials
will give us the scaled character $\chihat_{\boldptimesq}(k\; 1^{n-k})$.  In
fact, doing so produces polynomials in agreement with Stanley's data.  
We, therefore, can now use Kerov's polynomials to better understand Stanley's
character polynomials.  It is clear from \eqref{eq:stanphi} and
\eqref{eq:phipq} that $R_i(\boldptimesq)$ is a homogeneous polynomial of degree
$i$ in the $p$'s and $q$'s.  Therefore, since Kerov's polynomial
$\Sigma_k$ is graded with terms of weight $k+ 1 (\mathrm{mod}\; 2)$ (see
\citet[Proof of Theorem 1.1]{bi}) in the $R_i$'s, we see
that Stanley's character polynomials are also graded with terms of degree $k+1
(\mathrm{mod}\; 2)$.  We state this now as a proposition, for easy reference later.
\begin{proposition}\label{thm:kerovandstan}
	Terms of degree $i$ in $F_k(\boldp;\; \boldq)$ 
	are obtained from the terms of weight $i$ in Kerov's
	polynomials $\Sigma_k$ with the $R_i$'s evaluated at the shape
	$\boldptimesq$.
\end{proposition}
To further reinforce the idea that we are dealing with polynomials, and to make
convenient variable substitutions, we depart
from the notation used thus far.   We shall replace
$R_i(\boldptimesq)$ with $R_i(\boldpq)$ and $R_{\boldptimesq}(z)$ with
$R_\boldpq(z)$ to emphasize that these objects are polynomials in $p$'s and
$q$'s.  We do this analogously with 
$C(z)$,
$H_\boldptimesq(z)$ and $\phi_\boldptimesq(z)$;
that is, the series
$\phi_\boldpq(z)$,
will denote the series in \eqref{eq:phipq} and $H_\boldpq(z)$
will denote the series in \eqref{eq:stanphi}. 
We shall deal with the terms of different weights separately, starting
with the terms of highest degree, namely the terms of degree $k+1$.

\subsection{Terms of Degree $k+1$}\label{sec:kplusoneterms}

The expression for the terms of highest degree in Stanley's polynomials are given
implicitly by $G_{\boldpq}(z)$ in \eqref{eq:stantop}.  From Theorem \ref{thm:kerovtopterm} and Proposition
\ref{thm:kerovandstan}, we can obtain a similar
formula for the highest degree terms;  that is, the terms of highest degree in
$F_k(\boldp;\; \boldq)$, which have degree $k+1$,
is given by $R_{k+1}(\boldpq)$ and we see that the generating series for the terms of
highest degree is
\begin{equation}\label{eq:stantopkerov}
	R_{\boldpq}(z) = \f{z}{\left(\displaystyle \f{\displaystyle z \prod_{i=1}^m \left(1 - \left(q_i - \sum_{j=1}^i
	p_j\right)z\right)}{\displaystyle \left(1 + \sum_{j=1}^m p_j z\right) \prod_{i=1}^m \left(1-\left(q_i -
	\sum_{j=1}^{i-1} p_j\right) z\right)}\right)^\laginv}.
\end{equation}
Evidently, the two generating series $R_{\boldpq}(z)$ and $G_{\boldpq}(z)$
should be 
equal; after all they both generate the highest degree terms of $F_k(\boldp;\;
\boldq)$, although it is not obvious from \eqref{eq:stantop} and
\eqref{eq:stantopkerov} that this is the case.  It turns out
that $R_{\boldpq}(z)$ and $G_{\boldpq}(z)$ are \emph{almost} the same;  we state this
more precisely in the next proposition.
\begin{proposition}\label{thm:gkeqrk}
	The generating series $R_{\boldpq}(z)$ and $G_{\boldpq}(z)$ are
	identical except for the linear terms;  more precisely
	\begin{equation*}
		R_{\boldpq}(z) = G_{\boldpq}(z) - \sum_{i=1}^m p_i z.
	\end{equation*}
\end{proposition}
\begin{proof}
	From Lagrange inversion, it suffices to show that $R_{\boldpq}(z) +
	\sum_{i=1}^m p_i z$ satisfies the same equation as $G_{\boldpq}(z)$.
	In this proof, we denote $R_{\boldpq}(z)$ and
	$G_{\boldpq}(z)$ by $R$ and $G$, respectively, for convenience.  From
	\eqref{eq:stantopkerov} we have
	\begin{equation*}
		\f{z}{R} = \left(\displaystyle \f{\displaystyle z
		\prod_{i=1}^m \left(1 - \left(q_i - \sum_{j=1}^i
		p_j\right)z\right)}{\displaystyle \left(1 + \sum_{j=1}^m p_j
		z\right) \prod_{i=1}^m \left(1-\left(q_i - \sum_{j=1}^{i-1}
		p_j\right) z\right)}\right)^\laginv.
	\end{equation*}
	By the definition of compositional inverse we have, from the last
	expression,
	\begin{align*} 
		z &=
		\f{\displaystyle z \prod_{i=1}^m \left(R - \left(q_i -
		\sum_{j=1}^i p_j\right)z\right)}{\displaystyle \left(R +
		\sum_{j=1}^m p_j z\right) \prod_{i=1}^m \left(R-\left(q_i -
		\sum_{j=1}^{i-1} p_j\right) z\right)}\\ 
		&=
		\f{\displaystyle z \prod_{i=1}^m \left( \left(R + \sum_{j=1}^m
		p_j z \right) - \left(q_i + \sum_{j=i+1}^m p_j\right)
		z\right)}{\displaystyle \left(R + \sum_{j=1}^m p_j z\right)
		\prod_{i=1}^m \left(\left(R + \sum_{j=1}^m p_j z \right)
		-\left(q_i + \sum_{j=i}^{m} p_j\right) z\right)}\\
		&=
		\f{\displaystyle \f{z}{\left(R + \sum_{j=1}^m p_j z \right)}
		\prod_{i=1}^m \left( 1 - \left(q_i + \sum_{j=i+1}^m p_j\right)
		\f{z}{\left(R + \sum_{j=1}^m p_j z
		\right)}\right)}{\displaystyle 
		\prod_{i=1}^m \left(1 
		- \left(q_i + \sum_{j=i}^{m} p_j\right) \f{z}{\left(R +
		\sum_{j=1}^m p_j z \right)}\right)}.
	\end{align*}  
	Again, from the definition of compositional inverse, we conclude that
\begin{equation*}
	\f{z}{\left(R + \sum_{j=1}^m p_j z \right)} =
		\left( \f{\displaystyle z
		\prod_{i=1}^m \left( 1 - \left(q_i + \sum_{j=i+1}^m p_j\right)
		z\right)}{\displaystyle 
		\prod_{i=1}^m \left(1 
		- \left(q_i + \sum_{j=i}^{m} p_j\right) z \right)}
		\right)^\laginv.
\end{equation*}
	Comparing this expression with \eqref{eq:stantop}, the result follows.
\end{proof}
Indeed, using Lagrange inversion we see that the linear terms of $R_{\boldpq}(z)$ and
$G_{\boldpq}(z)$ are 0 and $\sum_{i = 0}^m p_i z$, respectively.
Furthermore,  note that although $R_{\boldpq}(z)$ and $G_{\boldpq}(z)$ differ in the linear term, this has no
effect on either Kerov's or Stanley's polynomials since $R_1(\boldpq)$ does not  
appear in Kerov's polynomials, as one can see in \eqref{eq:rexp} (in general,
this fact follows from the combinatorial argument given in \citet[Proof of
Theorem 1.1]{bi}).

Through Lagrange inversion, we see that the $R_i$ are written in terms of the
series $\phi_\boldpq$ given in \eqref{eq:altrk}.  We
use the notation $\phi_{\boldpminq},
R_k(\boldpminq)$ and $G_k(\boldpminq)$ to denote that we are substituting $-q_i$
for $q_i$ for all $i$ in these series.  We have the following compact
expression for the series $\phi_{\boldpminq}(-z)$.
\begin{theorem}\label{thm:phiform}
	For $\boldp = \seq{p}{m}$ and $\boldq = \seq{q}{m}$, we have
	\begin{equation*}
		\phi_{\boldpminq}(-z) = \prod_{i=1}^m \left(1
		+ \f{p_i q_i z^2}{(1 - r_{i-1} z) \left(1 - (q_i +
		r_i)z\right)} \right).
	\end{equation*}
	where $r_i = \sum_{j=1}^i p_j$.
\end{theorem}
\begin{proof}
	We have, from \eqref{eq:phipq},
	\begin{equation*}
		\phi_{\boldpminq}(-z) = \f{\displaystyle \left(1 - r_m z \right) \prod_{i=1}^m
		\left(1-\left(q_i + r_{i-1}
		\right)z\right)}{\displaystyle \prod_{i=1}^m
	\left(1 - \left(q_i + r_i \right)z\right)}.
	\end{equation*}
	Now set $A_n(z) = 1 - r_n z, F_0 = 1$ and
	\begin{equation}\label{eq:recurphi}
		F_n(z) = A_n(z) \f{\displaystyle \prod_{i=1}^n \left(1 -
		\left(q_i + r_{i-1}\right) z \right)}{\displaystyle \prod_{i=1}^n \left(1 -
		\left(q_i + r_i \right) z \right)}.
	\end{equation}
	Note that $\phi_{\boldpminq}(-z) = F_m(z)$.  Then,
	\begin{align}
		F_n(z) &= \f{F_{n-1}(z)}{A_{n-1}(z)} \f{1 - \left( q_n +
		r_{n-1} \right) z}{1 - \left( q_n + r_n \right)z} 
		A_n(z)\nonumber\\
		&= \f{F_{n-1}(z)}{A_{n-1}(z)} \f{\displaystyle A_{n-1}(z)
		\left(1 - \f{q_n z}{A_{n-1}(z)}\right)}{\displaystyle A_{n-1}(z)
		\left(1 - \f{\left( q_n + p_n \right)z}{A_{n-1}(z)}\right)}
		A_{n-1}(z) \left(1 - \f{p_n z}{A_{n-1}(z)}\right)\nonumber\\
		&= F_{n-1}(z) \f{\displaystyle 1 - \f{(q_n + p_n)z}{A_{n-1}(z)}
		+ \f{p_n q_n z^2}{A_{n-1}^2(z)}}{\displaystyle 1 - \f{(q_n +
		p_n)z}{A_{n-1}(z)}}\nonumber\\
		&= F_{n-1}(z) \left(1
		+ \f{p_n q_n z^2}{A_{n-1}^2(z) \left(1 - \f{(q_n +
		p_n)z}{A_{n-1}(z)}\right)} \right)\nonumber\\
		&= F_{n-1}(z) \left(1
		+ \f{p_n q_n z^2}{A_{n-1}(z) \left(1 - (q_n +
		r_n)z\right)} \right)\nonumber\\
		&= F_{n-1}(z) \left(1
		+ \f{p_n q_n z^2}{(1 - r_{n-1} z) \left(1 - (q_n +
		r_n)z\right)} \right).\label{eq:Falmost}
	\end{align}
	Therefore, from \eqref{eq:Falmost} we have
	\begin{align*}
		\phi_{\boldpminq}(-z) &= F_m(z)\\
		&= \f{F_m(z)}{F_0(z)}\\
		&= \f{F_m(z)}{F_{m-1}(z)} \cdot \f{F_{m-1}(z)}{F_{m-2}(z)}
		\cdots \f{F_1(z)}{F_{0}(z)}\\
		&= \prod_{i=1}^m \left(1
		+ \f{p_i q_i z^2}{(1 - r_{i-1} z) \left(1 - (q_i +
		r_i)z\right)} \right).
	\end{align*}
\end{proof}
\begin{cor}\label{thm:phipos}
	$\phi_{\boldpminq}(-z)$ is $\boldp,
	\boldq$-positive.
\end{cor}
\begin{proof}
	Each multiplicand in Theorem \ref{thm:phiform} is $\boldp,
	\boldq$-positive, making the product $\boldp, \boldq$-positive. 
\end{proof}
\begin{cor}\label{thm:rispos}
	For all $k \geq 1$, the series in $p$'s and $q$'s
	$(-1)^kR_{k+1}(\boldpminq)$ and\linebreak $(-1)^kG_k(\boldpminq)$ are $\boldp,
	\boldq$-positive.
	That is, the terms of highest degree in
	$(-1)^kF_k(\boldp;\; -\boldq)$ all have positive coefficients.
\end{cor}
\begin{proof}
	The series $(-1)^kG_k(\boldpminq)$ consists of by definition the terms of highest
	degree in $(-1)^kF_k(\boldp;\; -\boldq)$, and by Proposition
	\ref{thm:gkeqrk}, $(-1)^kG_k(\boldpminq) = (-1)^kR_{k+1}(\boldpminq)$
	are equal for all $k \geq 1$.  Thus, it suffices to show that
	$(-1)^kR_{k+1}(\boldpminq)$ is $\boldp, \boldq$-positive for all $k \geq
	1$.

	By \eqref{eq:altrk} we have
	\begin{align*}
		(-1)^k R_{k+1}(\boldpminq) &= (-1)^k \left(-\f{1}{k} [y^{k+1}]\;
		\phi_\boldpminq^k(y) \right)\\
		&= \f{1}{k} [(-y)^{k+1}] \phi_\boldpminq^k(y)\\
		&= \f{1}{k} [y^{k+1}] \phi_\boldpminq^k(-y),
	\end{align*} 
	and the result follows.
\end{proof}

\subsection{Terms of Degree $k-1$, $k-3$ and a General Connection Between
Kerov's Polynomials and Stanley's Polynomials}\label{sec:kminoneterms}

In this section we deal with terms of degree $k-1$ and $k-3$ in Stanley's
polynomials.  We note that in \citet{stan:7} there are no results concerning terms 
not of highest degree; Stanley comments only on the series
$G_{\boldpminq}(z)$, the terms of highest degree in $k+1$.  Moreover, we note the
complication that $(-1)^k \Sigma_k$ has some negative terms when one
evaluates the $R_i$ in terms of the shape $\boldpq$ and substitutes $-q_i$ for
all the $q_i$.  
More precisely, consider, for example, $\Sigma_5$ given in \eqref{eq:rexp}.  We see from the comments at the beginning of Section
\ref{sec:kerovstan} that
\begin{align*}
	(-1)^5 F_5(\boldpminq) &= (-1)^5 \Sigma_5(\boldpq)|_{\boldq \rightarrow
	-\boldq}\\
	&= (-1)^5 \left(R_6(\boldpminq) + 15R_4(\boldpminq) +
	5 R_2(\boldpminq)^2\right.\\
	& \hspace{3.0cm} \left. + 8 R_2(\boldpminq)\right)\\
	&= (-1)^5 R_6(\boldpminq) + 15 (-1)^3 R_4(\boldpminq)\\
	& \hspace{3.0cm} - 5 ( (-1)R_2(\boldpminq))^2 + 8 (-1)R_2(\boldpminq).
\end{align*}
Note that all terms are $\boldp, \boldq$-positive except for the term $-5 (
(-1)R_2(\boldpminq))^2$.  Thus, $\boldp, \boldq$-positivity would not
immediately follow from R-positivity of Kerov's polynomials.  For the terms of
degree $k-1$ and $k-3$, however, we can use Theorems \ref{thm:sigk2} and
\ref{thm:cpos}.  We begin with the following theorem.

\begin{theorem}\label{thm:stansecterms}
	For $k \geq 3$, the terms of degree $k-1$ in $F_k(\boldp;\; \boldq)$ are given by
	\begin{equation*}
		-\f{k (k+1)}{24} [y^{k-3}] \phi_\boldpminq^{\prime \prime}(y)
		\phi_\boldpminq^{k-1}(y).
	\end{equation*}
\end{theorem}
\begin{proof}
	From Proposition \ref{thm:kerovandstan} and Theorem \ref{thm:sigk2}, the terms of degree $k-1$ in $F_k(\boldp;\;
	\boldq)$ are given by
	\begin{equation*}
		\f{1}{4}{k+1 \choose 3} C_{k-1}(\boldpq).
	\end{equation*}
	Setting $w = z/R_{\boldpq}(z)$ then from \eqref{eq:recr} we
	have
	\begin{equation*}
		z = w R_{\boldpq}(z), \;\;\;\;\; w = z \phi_\boldpq(w),
	\end{equation*} 
	where $\phi_\boldpq(z)$ is given in \eqref{eq:phipq}.
	Further, from the definition of $C_{\boldpq}(z)$ we have
	\begin{equation*}
		C_{\boldpq}(z) = \f{1}{-z^2 \f{d}{dz} \f{1}{w}},
	\end{equation*}
	 Thus,
	\begin{equation*}
		z \f{d}{dz} w = \f{w}{1 - z 
		\phi_\boldpq^\prime(w)},
	\end{equation*}
	from which we obtain
	\begin{align*}
		C_{\boldpq}(z) &= \f{1}{-z^2 \f{d}{dz} \f{1}{w}}\\
		&= \f{1}{\f{z^2}{w^2} \f{d}{dz} w}\\
		&= \f{w}{z} (1 - z \phi_\boldpq^\prime(w))\\
		&= \phi_\boldpq(w) - w \phi_\boldpq^\prime(w).
	\end{align*}
	Therefore, for all $k \geq 2$, we have by Lagrange inversion that
	\begin{align*}
		[z^{k-1}]\; C_{\boldpq}(z) &= [z^{k-1}]\; \phi_\boldpq(w) -
		[z^{k-1}]\; w \phi_\boldpq^\prime(w)\\
		&= \f{1}{k-1} [y^{k-2}]\; \phi_\boldpq^\prime(y)
		\phi_\boldpq^{k-1}(y)\\
		&\phantom{doadeer} -
		\f{1}{k-1} [y^{k-2}]\; \left( \phi_\boldpq^\prime(y) + y
		\phi_\boldpq^{\prime \prime} (y) \right)
		\phi_\boldpq^{k-1}(y)\\
		&= - \f{1}{k-1} [y^{k-3}] \phi_\boldpq^{\prime \prime}(y)
		\phi_\boldpq^{k-1}(y),
	\end{align*}
	and the result follows.
\end{proof}
From Theorem \ref{thm:stansecterms} we obtain the following positivity result.
\begin{cor}\label{thm:corsectermspos}
	For $k \geq 3$, the terms of degree $k-1$ in $(-1)^k F_k(\boldp; -\boldq)$ are
	$\boldp, \boldq$-positive.
\end{cor}
\begin{proof}
	The terms of degree $k-1$ in $F_k(\boldp;\, \boldq)$ are
	given in Theorem \ref{thm:stansecterms}.  Therefore, the terms of degree
	$k-1$ in $(-1)^k F_k(\boldp;\; -\boldq)$ are
	\begin{align*}
		(-1)^k \f{1}{4}{k+1 \choose 3} C_{k-1}(\boldpminq) &= -\f{1}{4}
		{k+1 \choose 3} [z^{k-1}]\; C_{\boldpminq}(-z)\\
		&= \f{k (k+1)}{24} [y^{k-1}] (-y)^2 \f{d^2}{d(-y)^2} \left(
		\phi_{\boldpminq}(-y) \right)\\
		& \phantom{morebigtimegobbledeegook} \cdot \phi_{\boldpminq}^{k-1}(-y)\\
		&= \f{k (k+1)}{24} [y^{k-1}] y^2 \left(\f{d^2}{dy^2} (-1)^2
		\right) \left(
		\phi_{\boldpminq}(-y) \right)\\
		&\phantom{morebigtimegobbledeegook} \cdot \phi_{\boldpminq}^{k-1}(-y)\\
		&= \f{k (k+1)}{24} [y^{k-1}] y^2 \f{d^2}{dy^2} \left(
		\phi_{\boldpminq}(-y) \right) \phi_{\boldpminq}^{k-1}(-y).
	\end{align*}
	From Theorem \ref{thm:phipos} both $\phi_{\boldpminq}(-y)$ and, of
	course then, $\tf{d^2}{dy^2} \phi_{\boldpminq}(-y)$ are $\boldp,
	\boldq$- \linebreak positive, proving the result.
\end{proof}
The following theorem gives a general connection between Kerov's polynomials and
Stanley's polynomials.  
\begin{theorem}\label{thm:connectkerovstan}
	If Kerov's polynomials $\Sigma_k$ are C-positive then Stanley's
	polynomials
	\linebreak
	$(-1)^k F_k(\boldp;\; -\boldq)$ are $\boldp,
	\boldq$-positive.
\end{theorem}
\begin{proof}
	From Proposition \ref{thm:kerovandstan} the terms of degree $i$ in
	Stanley's polynomials are obtained from the terms of weight $i$ in
	Kerov's polynomials.  From Theorem \ref{thm:gencexp} the terms of
	degree $k+1 -2n$ in Stanley's polynomials are obtained from
\begin{equation*}
		\sum_{i_1, \dots, i_{2n-1} \geq 0 \atop i_1
		+ \cdots + i_{2n-1} = k+1 -2n} \gamma_{i_1,
		\dots, i_{2n-1}} C_{i_1}(\boldp;\; \boldq) \cdots
		C_{i_{2n-1}}(\boldp;\; \boldq).
\end{equation*}	
Thus, the terms of degree $k+1 -2n$ in $(-1)^k F_k(\boldp;\; -\boldq)$ are
given by
\begin{equation*}
		\sum_{i_1, \dots, i_{2n-1} \geq 0
		\atop i_1 + \cdots + i_{2n-1} = k+1 -2n} \gamma_{i_1, \dots,
		i_{2n-1}} \left((-1)^{i_1 -1}C_{i_1}(\boldp;\; -\boldq)\right)
		\cdots \left((-1)^{i_{2n-1} -1} C_{i_{2n-1}}(\boldp;\;
		-\boldq)\right).
\end{equation*}
From the proof of Corollary \ref{thm:corsectermspos}, each
$(-1)^{j-1}C_j(\boldp;\; -\boldq)$ is $\boldp, \boldq$-positive, and the result
follows.
\end{proof}
\begin{cor}\label{thm:stanthirdterms}
	For $k \geq 5$, the terms of  degree $k-3$ in $(-1)^k F(\boldp;
	-\boldq)$ are $\boldp, \boldq$-positive.
\end{cor}
\begin{proof}
	Follows directly from Theorems \ref{thm:cpos} and
	\ref{thm:connectkerovstan}.
\end{proof}

\section*{Acknowledgements}

I would like to thank Ian Goulden for suggesting this topic, in particular for
pointing out the paper \citet{stan:7}.  I would also like to thank John Irving
for very useful conversations concerning this material.\cite{ratt:2}
\bibliographystyle{plainnat}
\bibliography{bibdatabase}

\end{document}